\documentclass[letterpaper, 10 pt, conference]{ieeeconf}  

\IEEEoverridecommandlockouts                              

\usepackage{amsmath}
\usepackage{amssymb}
\usepackage{amsfonts}
\usepackage{mathrsfs}
\usepackage{graphicx}
\usepackage{overpic}

\usepackage{color}
\usepackage{soul}

\usepackage{theorem}

\newtheorem{lemma}{Lemma}
\newtheorem{remark}{Remark}
\newtheorem{assumption}{Assumption}
\newtheorem{proposition}{Proposition}

\newcommand{\vect}[1]{\ensuremath{\boldsymbol{\mathrm{#1}}}}

\definecolor{seb}{rgb}{0.8,1,1}

\definecolor{kai}{rgb}{0.61, 0.87, 1.0}

\definecolor{todo}{rgb}{0.68, 1.0, 0.18} 

\newcommand{\BeginProof}{\noindent \textbf{Proof.}$\;\;$}
\newcommand{\EndProof}{\hfill{$\blacksquare$}\vspace{0.3cm}}

\title{\LARGE \bf
	Solving Mission-Wide Chance-Constrained Optimal Control Using Dynamic Programming
}
\author{Kai Wang and S\'ebastien Gros 
	\thanks{The authors are with the Department of Engineering Cybernetics, Norwegian University of Science and Technology (NTNU), 7491, Trondheim, Norway. {\tt\small\{kai.wang, sebastien.gros\}@ntnu.no} }
}

\begin{document}

\maketitle
\thispagestyle{empty}
\pagestyle{empty}

%

\begin{abstract}

	This paper aims to provide a Dynamic Programming (DP) approach to solve the Mission-Wide Chance-Constrained Optimal Control Problems (MWCC-OCP). The mission-wide chance constraint guarantees that the probability that the entire state trajectory lies within a constraint/safe region is higher than a prescribed level, and is different from the stage-wise chance constraints imposed at individual time steps. The control objective is to find an optimal policy sequence that achieves both (i) satisfaction of a mission-wide chance constraint, and (ii) minimization of a cost function. By transforming the stage-wise chance-constrained problem into an unconstrained counterpart via Lagrangian method, standard DP can then be deployed. 
	Yet, for MWCC-OCP, this methods fails to apply, because the mission-wide chance constraint cannot be easily formulated using stage-wise chance constraints due to the time-correlation between the latter (individual states are coupled through the system dynamics). To fill this gap, firstly, we detail the conditions required for a classical DP solution to exist for this type of problem; secondly, we propose a DP solution to the MWCC-OCP through state augmentation by introducing an additional functional state variable.
\end{abstract}

\section{INTRODUCTION}\label{sec:intro}
In real-world control applications uncertainties, plant-model mismatch and exogenous disturbances are unavoidable. Some finite-horizon Markov control problems require that a risk criterion must be met over the entire planning horizon. Such a requirement could arise due to certain regulations, such as the safety-stock placement in supply chain management and the collision-avoidance in robot path planning. Performing optimal control of such system under hard (robust) constraints can yield very conservative control policies or even be infeasible. Another commonly used approach is to add a penalty term to the objective function. However, it is hard to decided how strict the penalty should be set; high penalties may result in conservative solutions, while low penalties may result in high risk.  It is therefore common to rather handle the risk using chance (probabilistic) constraints directly.

In the literature of optimal control, there are many forms of chance constraint. E.g., in path planning for vehicles in the presence of obstacles, the mission is supposed to plan an optimal trajectory for a vehicle within $N$ time stages. Suppose that there exist $m$ constraints for the vehicle at each state and time stage. \textit{Individual chance constraint} restricts the probability that the state violates a single constraint at a single time stage. \textit{Stage-wise chance constraint} restricts the probability that the state violates any individual constraints at a single time stage, which consists of $m$ individual chance constraints. Both of them in effect restrict at every time stage the probability that the vehicle collides with an obstacle. \textit{Mission-wide chance constraint} restricts the probability that the state sequence/trajectory spanning over all time steps violates any constraints among the total $mN$ constraints. In contrast to individual or stage-wise ones, a mission-wide chance constraint directly restricts the probability of collision on the overall driving mission. A mission-wide chance constraint is arguably more meaningful than stage-wise constraints in some specific control tasks. Indeed, the former directly handles the risk of running a mission \cite{Lew2019,Ono08,Ono2015}, while the latter does it very indirectly since satisfying a risk level in each time step may result in a poor risk level over the entire horizon, see Section II-A of \cite{Wang2021} for detailed explanations. However, individual or stage-wise chance constraints are easier to handle than mission-wide constraints. Indeed, the mission-wide chance constraints involve probabilities over entire state trajectory, yielding very large probability spaces. 

Notice that there is another widely used terminology, \textit{joint chance constraint}, referring to the stage-wise constraint in e.g. \cite{Mesbah2018}, as well as the mission-wide constraint in e.g. \cite{Lew2019}, and even the conjunction of the stage-wise constraints within only the prediction horizon that is typically shorter than the mission duration in some stochastic Model Predictive Control (MPC) literature, e.g., \cite{Paulson2020}. In this paper we do not adopt this terminology to reduce the risk of confusion. 

Optimization subject to the chance constraints was first proposed in the seminal work \cite{Charnes1959}. Chance-constrained optimization problems are typically intractable. Two main reasons behind this intractability are (i) the difficulty of checking the feasibility of a solution as it requires evaluating multivariate integrals, and (ii) the non-convex feasible region defined by a chance constraint.    

The current research on Chance-Constrained Optimal Control Problems (CC-OCP) is centered around tractable approximation approaches, such as Convex bounding approaches  \cite{Ono08,Van2016,Paulson2020} and sampling approaches \cite{Campi2006,Nemirovski2007}. The closed-loop controller, stochastic MPC, has been widely used in practice to implement the CC-OCP and has been intensively investigated among the MPC community, see \cite{Mesbah2016,Kouvaritakis2016} and references therein. It approximates the optimal control policy sequence and executes its control action online. To the best of our knowledge, however, there does not exist any stochastic MPC schemes that explicitly guarantee a rigorous satisfaction of the mission-wide chance constraint. \cite{Wang2021} provides a tentative stochastic MPC solution that ensures recursive feasibility in a specific sense, while the mission-wide chance constraint satisfaction is conserved.

Stochastic Dynamic Programming (DP) is a general framework for model-based sequential decision-making processes under uncertainty and provides a global optimal control policy sequence. Chance-constrained DP was first introduced in \cite{Askew1974} with application in operation of reservoir, and was extended in \cite{Sniedovich1975} by introducing additional system variables. In \cite{Askew1974,Sniedovich1975}, however, these approaches are heuristic and the chance constraints are limited to the form of cumulative stage-wise chance constraints. \cite{Rossman1977} provided a systematic way to cope with the formulation proposed in \cite{Askew1974} within the framework of Lagrangian duality theory, and also pointed out that there is no clear way to incorporate the mission-wide chance constraint because of the lose of additive structure in the constraints. This is due to the fact that there exists time-correlation between the stage-wise chance constraints such that the mission-wide chance constraint cannot be simply expressed as the summation of stage-wise chance constraints. The authors in \cite{Ono2012} proposed a DP method for MWCC-OCP, by first reformulating the mission-wide chance constraint via Boole's approximation, and then solving the resulting unconstrained Lagrangian function. This method was applied to the robotic space exploration mission in \cite{Ono2015}. When the optimization criteria is absent, in \cite{Pola2006}, a tailored DP approach was proposed to maximize the mission-wide probability of safety. Unfortunately, MWCC-OCP cannot necessarily be put in that simple form.

The above statement does not exclude the existence of a DP solution to the MWCC-OCP, but it indicates that classical DP framework cannot be simply applied. This paper entails in theory that there does exist an exact DP solution to the MWCC-OCP through proper state augmentation. This comes at the cost of the augmented state being in a functional space, hence making the resulting DP problem significantly more challenging to tackle.

The paper is structured as follows. In Section \ref{sec:statement} we formulate the MWCC-OCP and provide some preliminaries. Section \ref{sec:classicalDP} is devoted to develop a DP algorithm for the general risk-constrained dynamic optimization problem. In Section \ref{sec:augmentation0} we detail a specific DP algorithm to the MWCC-OCP through state augmentation. An one-dimensional lineal case study is given in Section \ref{sec:casestudy}. Finally, Section \ref{sec:conclusion} concludes the paper.

\section{Problem setup and preliminaries}\label{sec:statement}
We will consider a discrete-time Markov Decision Process (MDP) with
continuous state and action space. More specifically, we assume the state space $\mathcal S\subseteq\mathbb R^{n_s}$ and control/action space $\mathcal A\subseteq\mathbb R^{n_a}$. The transition density function  
\begin{equation}\label{eq:PDF}
	\rho[\,\vect s^{+}\,|\, \vect s,\vect a\,]
	:\mathcal S\times\mathcal S\times\mathcal A\to [0,\infty),
\end{equation}
provides a conditional probability of observing a transition from the state-action pair $\vect s\in\mathcal S,\,\vect a\in\mathcal A$ to a successor state $\vect s^+ \in\mathcal S$. In the control community, the system dynamics are often alternatively defined as $\vect s^{+} = \vect f(\vect s, \vect a, \vect w),$
where $\vect f:\mathcal S\times\mathcal{A}\times\mathcal W\to\mathcal{S}$ is possibly a vector-valued nonlinear function and the disturbance $\vect w$ is the realization of a random variable that takes values from set $\mathcal W \subseteq\mathbb{R}^{n_w}$. 
We will assume that the finite stage cost is given by $\ell(\vect s, \vect a):\mathcal{S}\times\mathcal{A}\to \mathbb R$, 
and that deterministic, Markovian policies are functions $\vect{\pi}_k(\vect s_k):\mathcal{S}\to\mathcal{A}$,
which specify an action $\vect a_k=\vect{\pi}_k(\vect s_k)$
when the system visits state $\vect s_k$ at time step $k$. 

Let a nonzero $N\in \mathbb N$ denote the mission duration. To accomplish a mission,  a decision maker or controller has to make actions/decisions that minimize the following cost
\begin{equation}
	\mathbb E \left[ \sum_{k=0}^{N-1}\ell(\vect s_k,\vect a_k) + \ell_N(\vect s_N)\,\middle|\,\vect s_0,\vect a_k=\vect{\pi}_k(\vect s_k) \right]\,.
\end{equation}  
Here, $\mathbb E$ is the expected value operator applying to the space of all possible trajectories $\vect s_{0,\ldots,N}$ of the states in closed-loop with policy sequence $\vect\pi := \{\vect\pi_0,\ldots,\vect\pi_{N-1}\}$ and a fixed initial state $\vect s_0\in\mathbb S$.
Function $\ell_N$ is a finite terminal cost incurred at the end of the mission. We set $\vect\pi^k := \{\vect\pi_k,\ldots,\vect\pi_{N-1}\}$
for any $k=1,\ldots,N-1$. 

We consider a mission to be safe if the closed-loop trajectory of state
lies in a constraint set $\mathbb S\subseteq\mathcal{S}$, i.e., $\vect s_{1,\ldots,N}\in\mathbb{S}.$ 
However, in many practical cases, satisfying this safety requirement may be impossible especially when the disturbances of the system are unbounded. Even in the case where the system states are finitely supported, it can still be infeasible, and yield very conservative policies. Alternatively, we seek to guarantee a probabilistic safety on the state sequence spanning over the whole mission, i.e., a mission-wide chance constraint (MWCC):
\begin{equation}\label{eq:MWCC}
	\mathbb{P}\left[\,\vect s_{1,\ldots,N}\in\mathbb S\,\middle|\,\vect s_0,\vect\pi \right] \geq 1-\varepsilon\,,
\end{equation} 
where  $\mathbb P[\cdot]$ denotes the probability operator, and $\varepsilon\in[0,1]$ is a predefined \textit{risk bound} indicating that the probability that the states stay within a safe set over the mission duration or planning horizon is at least $1-\varepsilon$.  We call the left hand term of \eqref{eq:MWCC}, i.e., $\mathbb{P}\left[\, \vect s_{1,\ldots,N}\in\mathbb S\,\middle|\,\vect s_0,\vect\pi \right],$ Mission-Wide Probability of Safety (MWPS).

The resulting MWCC-OCP can be formally formulated as
\begin{subequations}\label{P:original}
	\begin{align} 
		J^{*}(\vect s_0)=\min_{\vect\pi}&\,\,\, \mathbb E \left[ \sum_{k=0}^{N-1}\ell(\vect s_k,\vect a_k) + \ell_N(\vect s_N)\,\middle|\,\vect s_0,\vect\pi  \right]\label{P:original_1}\\ 
		\mathrm{s.t.}&\,\,\,\mathbb{P}\left[\, \vect s_{1,\ldots,N}\in\mathbb S\,\middle|\,\vect s_0,\vect\pi \right] \geq 1-\varepsilon\,,\label{P:original_2}
	\end{align}
\end{subequations}
for each feasible initial state $\vect s_0\in\mathbb S$. Here we assume that $\vect{s}_0$ lies in the set from which there exists police sequence $\vect\pi$ such that \eqref{P:original_2} is feasible, even through such set is not easy to be calculated \cite{Vinod2021}. We will call $J^*$ the optimal cost function that assigns the optimal cost $J^{*}(\vect s_0)$ to each feasible initial state $\vect s_0$.


To solve problem \eqref{P:original}, apart from the approximation methods mentioned in Section \ref{sec:intro}, one would naturally think of DP because it is arguably the most general approach to sequential decision-making problem when a model is at hand. Even when it is computationally expensive, DP may still serve as the basis for many practical approaches. Since the above mission-wide chance constraint are acting on the whole Markov Chain $\vect s_{1,\ldots,N}$, existing methods for constrained DP \cite{Bertsekas05} and constrained MDPs \cite{Altman1999} cannot be used. This is due to the fact that these methods require that the state constraints are expressed in additive or independent form. As we will show in this paper, problem \eqref{P:original} presents interesting features.

In the following section, we integrate the constraint of problem \eqref{P:original} into its objective function via a penalty function, and then investigate conditions under which a DP scheme can be deployed.

\section{MWPS-constrained problem with DP}\label{sec:classicalDP} 

Let $\mathbf{1}_{\mathbb S}: \mathcal{S}\to\{0,\,1\}$ denote the indicator function of set $\mathbb S\subseteq\mathcal S;$ $\mathbf{1}_{\mathbb S}(\vect s) =1$, if $\vect s\in\mathbb{S}$, and $0$ if $\vect s\notin\mathbb{S}$.
Let us define the set of functions $V^{\vect\pi^k}_k:\mathbb{S}\to[0,1],\,\forall k=0,1\ldots,N-1$ as follows:
\begin{align*}
	&V^{\vect\pi^N}_N (\vect s)=V_N (\vect s) = \mathbf{1}_{\mathbb S}(\vect s)\,,\\
	&V^{\vect\pi^k}_k (\vect s)=\mathbb{P}\left[\,\vect s_{k+1,\ldots,N}\in\mathbb S \,\middle|\,\vect s_k=\vect s, \vect \pi^k\,\right].
\end{align*} 
In this way, given a state $\vect s_0$, the associated MWPS can be denoted by $V^{\vect\pi}_0(\vect s_0)=V^{\vect\pi^0}_0(\vect s_0)=\mathbb{P}\left[\, \vect s_{1,\ldots,N}\in\mathbb S \,\middle|\,\vect s_0,\vect\pi \right].$

\begin{lemma}\label{lem:BackRecur}
	Function $V^{\vect\pi^k}_k$ can be computed by
	the backward recursion:
	\begin{align}\label{eq:BackRecur}
		V^{\vect\pi^k}_{k} (\vect s_{k})&= \int_{\mathbb S}V^{\vect\pi^{k+1}}_{k+1}(\vect s_{k+1})\cdot \rho[\vect s_{k+1}\mid\vect s_k, \vect \pi_k(\vect s_k)]\,\mathrm{d}\vect s_{k+1}\nonumber\\
		&:=\mathbb E_{\vect s_{k+1}}\left[V^{\vect\pi^{k+1}}_{k+1}(\vect s_{k+1})\middle | \vect s_k, \vect\pi_k(\vect s_k) \right]
	\end{align}
	initialized with the boundary condition $V_N (\vect s_N) = \mathbf{1}_{\mathbb S}(\vect s_N),$
	where operation $\mathbb E_{\vect s_{k+1}}$ indicates that the expectation is taken with respect to the probability distribution of $\vect s_{k+1}$ that remains in set $\mathbb S$.
\end{lemma}
\BeginProof See \cite[Lemma 1]{Abate2008}\EndProof 

As what is typically done in constrained MDPs, penalizing the risk of mission failures with a
suitably chosen cost function can in principle guarantee that the MDPs yield a
policy which tends to not violate the constraints or even not violate the constraints at all when exact penalty is used.

let us consider an optimization problem of the form: 
\begin{align}
	&\tilde{J}^*(\vect s_0)=\min_{\vect\pi}\,\,\mathbb E \left[ \sum_{k=0}^{N-1}\ell(\vect s_k,\vect a_k) + \ell_N(\vect s_N)\,\middle|\,\vect s_0,\vect\pi  \right]\nonumber\\
	&\hspace{5.3cm}+\zeta\left( V^{\vect\pi}_0(\vect s_0)\right), \label{eq:OptProb}
\end{align} 
where $\zeta:[0,1]\to[-\infty,+\infty]$ is a penalty function mapping the MWPS into a scalar that possibly takes values in the extended real line.

\begin{assumption}\label{ass:commutes}
	The function $\zeta$ commutes with the expectation operator $\mathbb E_{\vect s_{k+1}}$ for all $k=0,1,\ldots,N-1$, i.e., 
	\begin{align}\label{eq:CommuteAssum}
		\zeta\left( V^{\vect\pi^{k}}_{k}(\vect s_{k}) \right) =&\,\zeta\left( \mathbb E_{\vect s_{k+1}}\left[V^{\vect\pi^{k+1}}_{k+1} (\vect s_{k+1})  \middle | \vect s_k, \vect\pi_k(\vect s_k)\right]\right) \nonumber \\	
		=&\,\mathbb E_{\vect s_{k+1}}\left[ \zeta\left(V^{\vect\pi^{k+1}}_{k+1} (\vect s_{k+1})\right)\middle | \vect s_k, \vect\pi_k(\vect s_k)\right]
		\end{align}	
\end{assumption}
\begin{remark}
		It is obvious to verify that Assumption \ref{ass:commutes} holds if function $\zeta$ is affine and may only apply for this case. The commutation property in this assumption is essentially what we need to derive the following proposition. 
\end{remark}

\begin{proposition}\label{prop:ClassicalDP}
	Let Assumption \ref{ass:commutes} be satisfied. Then \eqref{eq:OptProb} can be solved via the following  DP recursion on the state-space:
	\begin{align}
		&\tilde{J}_N(\vect s_N) = \ell_N(\vect s_N)+\zeta\left(V_N(\vect s_N)\right) \,,\nonumber\\
		&\tilde{J}_k(\vect s_k)\hspace{-0.07cm} = \hspace{-0.07cm}\min_{\vect a_k\in\mathcal A}\ell(\vect s_k,\vect a_k)\hspace{-0.05cm}+\hspace{-0.1cm}\int_{\mathbb S}\hspace{-0.05cm} \tilde{J}_{k+1}(\vect s_{k+1})\rho[\vect s_{k+1}|\vect s_k,\vect a_k]\mathrm{d} \vect s_{k+1}\nonumber\\
		&\hspace{4cm}k=N-1,N-2,\ldots,0.\label{eq:NormDP_2}
	\end{align}
For every initial state $\vect s_0$, the optimal cost $\tilde{J}^*(\vect s_0)$ of problem \eqref{eq:OptProb} is equal to $\tilde{J}_0(\vect s_0)$, given by the last step (backward in time) of the above recursion. Furthermore, if $\vect\pi_k^*(\vect s_k)=\vect a_k^*$ minimizes the right hand side of \eqref{eq:NormDP_2} for each $\vect s_k$ and $k$, then the policy $\vect\pi^*=\{\vect\pi^*_0,\ldots,\vect\pi^*_{N-1}\}$ is optimal for problem \eqref{eq:OptProb}.
\end{proposition}
\BeginProof
For $k=N-1,N-2,\ldots,0$, let $\tilde{J}^*_k(\vect s_k)$ be the optimal cost for the $(N-k)$-stage problem associated to \eqref{eq:OptProb} that starts at state $\vect s_k$ and time $k$, and ends at time $N$,
\begin{align*}
	&\tilde{J}^*_k(\vect s_k) = \min_{\vect \pi^k}\,\mathbb E \left[ \sum_{i=k}^{N-1}\ell(\vect s_i,\vect a_i) + \ell_N(\vect s_N)\,\middle|\,\vect s_k,\vect\pi^k  \right]\\
	&\hspace{5.0cm}+\zeta\left(V^{\vect\pi^k}_k (\vect s_k) \right).
\end{align*}  
For $k = N$, we define: 
$\tilde{J}^*_N(\vect s_N) = \ell_N(\vect s_N)+\zeta\left(V_N(\vect s_N)\right).$
We will show by induction that the functions $\tilde{J}^*_k$ are equal to the functions $\tilde{J}_k$ generated by the DP recursion \eqref{eq:NormDP_2}, such that for $k=0$ the desired result will be obtained. 

By definition, we have $\tilde{J}^*_N = \tilde{J}_N = \ell_N+\zeta.$ Assume that for some $k+1$ and all $\vect s_{k+1}$, we have $\tilde{J}^*_{k+1}(\vect s_{k+1})=\tilde{J}_{k+1}(\vect s_{k+1})$. Then, since $\vect\pi^k=(\vect\pi_k,\vect\pi^{k+1})$, we have for all $\vect s_k\in\mathbb{S}$, (these developments are further explained hereafter)
\begin{align*}
	&\tilde{J}^*_k(\vect s_k) = \min_{\vect\pi^k}\,\mathbb E \left[ \sum_{i=k}^{N-1}\ell(\vect s_i,\vect a_i) + \ell_N(\vect s_N)\,\middle|\,\vect s_k,\vect\pi^k  \right]\\
	&\hspace{5cm}+\zeta\left(V^{\vect\pi^k}_k (\vect s_k) \right)\\
	&=\min_{\vect\pi_k}\,\ell(\vect s_k,\vect \pi_k(\vect s_k))\\
	&\hspace{0.0cm}+\mathbb{E}_{\vect s_{k+1}}\hspace{-0.1cm}\left[\min_{\vect\pi^{k+1}}\mathbb E \hspace{-0.05cm}\left[ \sum_{i=k+1}^{N-1}\ell(\vect s_i,\vect a_i) + \ell_N(\vect s_N)\,\middle|\,\vect s_{k+1},\vect\pi^{k+1}  \right]\right]\\
	&+\min_{\vect\pi^{k+1}} \zeta\left( \mathbb E_{\vect s_{k+1}}\left[V^{\vect\pi^{k+1}}_{k+1} (\vect s_{k+1}) \middle | \vect s_k, \vect\pi_k(\vect s_k) \right]\right).
\end{align*}
Based on equality \eqref{eq:CommuteAssum} in Assumption \ref{ass:commutes}, by commuting $\mathbb{E}_{\vect s_{k+1}}$ and $\zeta$, we can further observe that the above equation is equal to  
\begin{align*}
	&\tilde{J}^*_k(\vect s_k) =\min_{\vect\pi_k}\,\ell(\vect s_k,\vect \pi_k(\vect s_k))\\
	&\hspace{0.0cm}+\mathbb{E}_{\vect s_{k+1}}\hspace{-0.1cm}\left[\min_{\vect\pi^{k+1}}\mathbb E \hspace{-0.05cm}\left[ \sum_{i=k+1}^{N-1}\ell(\vect s_i,\vect a_i) + \ell_N(\vect s_N)\,\middle|\,\vect s_{k+1},\vect\pi^{k+1}  \right]\right]\\
	&+\min_{\vect\pi^{k+1}} \mathbb E_{\vect s_{k+1}}\left[\zeta\left( V^{\vect\pi^{k+1}}_{k+1} (\vect s_{k+1}) \right)\middle | \vect s_k, \vect\pi_k(\vect s_k)\right]
\end{align*}
Furthermore, we observe that 
\begin{align*}
	&\tilde{J}^*_k(\vect s_k) =\min_{\vect\pi_k}\,\ell(\vect s_k,\vect \pi_k(\vect s_k))\\
	&+\mathbb{E}_{\vect s_{k+1}}\hspace{-0.1cm}\left[\min_{\vect\pi^{k+1}}\mathbb E \hspace{-0.05cm}\left[ \sum_{i=k+1}^{N-1}\ell(\vect s_i,\vect a_i) + \ell_N(\vect s_N)\,\middle|\,\vect s_{k+1},\vect\pi^{k+1}  \right]\right]\\
	&+ \mathbb E_{\vect s_{k+1}}\left[\min_{\vect\pi^{k+1}}\zeta\left( V^{\vect\pi^{k+1}}_{k+1} (\vect s_{k+1}) \right) \middle | \vect s_k, \vect\pi_k(\vect s_k)\right]\\[0.3cm]
	&=\min_{\vect\pi_k}\,\ell(\vect s_k,\vect \pi_k(\vect s_k))\\
	&\hspace{0.2cm}+\mathbb{E}_{\vect s_{k+1}}\hspace{-0.1cm}\left[\min_{\vect\pi^{k+1}}\mathbb E \hspace{-0.05cm}\left[ \sum_{i=k+1}^{N-1}\ell(\vect s_i,\vect a_i) + \ell_N(\vect s_N)\,\middle|\,\vect s_{k+1},\vect\pi^{k+1}  \right]\right.\\
	&\hspace{2.3cm}\left.+\zeta\left( V^{\vect\pi^{k+1}}_{k+1} (\vect s_{k+1}) \right)\,\middle |\, \vect s_k, \vect\pi_k(\vect s_k)\right] .
\end{align*}
In the first equation above we moved the second $\min_{\vect\pi^{k+1}}$ inside the brackets expression $\mathbb E_{\vect s_{k+1}}[\,\cdot\,]$. In the last equation above, we integrated the two expressions of $\mathbb E_{\vect s_{k+1}}[\,\cdot\,]$ into one. Finally, by using the definition of $\tilde{J}^*_{k+1}(\vect s_{k+1})$, we can further have 
\begin{align*}
	&\tilde{J}^*_k(\vect s_k)\\
	&= \min_{\vect\pi_k}\,\ell(\vect s_k,\vect \pi_k(\vect s_k)) + \mathbb E_{\vect s_{k+1}}\left[\tilde{J}^*_{k+1}(\vect s_{k+1})\middle | \vect s_k, \vect\pi_k(\vect s_k)  \right] \\	
	&= \min_{\vect\pi_k}\,\ell(\vect s_k,\vect \pi_k(\vect s_k)) + \mathbb E_{\vect s_{k+1}}\left[\tilde{J}_{k+1}(\vect s_{k+1}) \middle | \vect s_k, \vect\pi_k(\vect s_k) \right] \\
	&= \min_{\vect a_k\in\mathcal A}\,\ell(\vect s_k,\vect a_k) + \int_{\mathbb S}\tilde{J}_{k+1}(\vect s_k) \rho\left[\vect s_{k+1}\,\middle|\,\vect s_k, \vect a_k\right]\mathrm{d}\vect s_{k+1}\\
	&= \tilde{J}_k(\vect s_k)\,.
\end{align*}
In the second equation, we used the induction hypothesis. In the third equation, we converted the minimization over $\vect\pi_k$ to a minimization over $\vect a_k$, using the fact that for any function $\kappa$ of $\vect s$ and $\vect a$, we have 
\[
\min_{\vect\mu\in\Pi}\kappa(\vect s, \vect\mu(\vect s)) = \min_{\vect a\in\mathcal A}\kappa(\vect s, \vect a)\, 
\]
where $\Pi$ is the set of all functions $\vect\mu(\vect s)$ such that $\vect\mu(\vect s)\in\mathcal A$ for all $\vect s$.    
\EndProof

\begin{remark}\label{remark:linear}
	If $\zeta$ is an affine function of the form $\zeta(x) = \lambda x + \Delta,\quad\forall\,x\in[0,1] , $ 
	where $\lambda, \Delta\in\mathbb R$ are constants, it is easy to verify that this class of functions satisfy equality \eqref{eq:CommuteAssum} in Assumption \ref{ass:commutes}. Indeed, we observe that 
	 \begin{align*}
	 	\zeta\left( V^{\vect\pi^{k}}_{k} (\vect s_{k})\right) =
	 	&\,\zeta\left( \mathbb E_{\vect s_{k+1}}\left[V^{\vect\pi^{k+1}}_{k+1} (\vect s_{k+1})  \middle | \vect s_k, \vect\pi_k(\vect s_k)\right]\right)\\
	 	= &\lambda\mathbb E_{\vect s_{k+1}}\left[V^{\vect\pi^{k+1}}_{k+1} (\vect s_{k+1})  \middle | \vect s_k, \vect\pi_k(\vect s_k)\right] +\Delta \\
	 	=&\mathbb E_{\vect s_{k+1}}\left[ \lambda V^{\vect\pi^{k+1}}_{k+1} (\vect s_{k+1}) + \Delta \middle | \vect s_k, \vect\pi_k(\vect s_k)\right]\\
	 	=&\mathbb E_{\vect s_{k+1}}\left[ \zeta \left( V^{\vect\pi^{k+1}}_{k+1} (\vect s_{k+1}) \right)\middle | \vect s_k, \vect\pi_k(\vect s_k)\right].
	 \end{align*}
\end{remark}
\begin{remark}
		In the literature on stochastic reachability analysis, see e.g.\cite{Abate2008}, if the first term
		in the cost functions of \eqref{eq:OptProb} is absent, a DP approach is proposed in \cite{Abate2008} to minimize the risk only, that is, $\min_{\vect\pi}\,\,\zeta\left( V^{\vect\pi}_0(\vect s_0)\right) :=   1-V^{\vect\pi}_0(\vect s_0)$,
		which is a special case of problem \eqref{eq:OptProb}. Here, the parameters of $\zeta$ are given by $\lambda=-1$ and $\Delta=1$. 
\end{remark}

To further make problem \eqref{eq:OptProb} equivalent to the original problem \eqref{P:original}, that is, $\tilde{J}^*(\vect s_0) = \tilde{J}(\vect s_0)$ for all feasible $\vect s_0$, we need to specify a concrete function $\zeta$. 

\begin{proposition}\label{pro:equiv}
	Suppose that function $\zeta$ is of the form:
	\begin{align}\label{eq:zeta}
		\zeta\left(x \right) = 
		\begin{cases}
			0,        &\text{if}\,\, x \geq 1-\varepsilon \\
			\infty,   &\text{otherwise}\,.
		\end{cases}
	\end{align} We observe that problem \eqref{P:original} and problem \eqref{eq:OptProb} are equivalent.
\end{proposition}
\BeginProof 
This proof is straightforward. Suppose that there is a set of policy sequences such that the MWCC \eqref{eq:MWCC} is satisfied. An optimal solution to problem \eqref{eq:OptProb} must fall in this set, because any policy sequence outside this set will yield an infinite penalty and hence cannot be optimal. Moreover, all the polices in this set result in zero penalty. Thus, the proof is concluded.   
\EndProof

However, Assumption \ref{ass:commutes} is not satisfied by function $\zeta$ defined in \eqref{eq:zeta} because in this case $\zeta$ does not commute with the expectation operator $\mathbb{E}_{\vect s_{k+1}}[\,\cdot\,]$. This means that in the case of function $\zeta$ given by \eqref{eq:zeta}, Proposition \ref{prop:ClassicalDP} is no longer applicable.

Note that the above observation implies that problem \eqref{eq:OptProb} with $\zeta$ defined in \eqref{eq:zeta}, which is equivalent to the original problem \eqref{P:original}, needs to be further reformulated into the basic problem format where DP algorithm can be deployed. This will be the subject of the following section.

\section{Functional State Augmentation}\label{sec:augmentation0}

The observations of the previous section show that adding penalty function associated to the MWPS into the objective function is of limited use if that penalty is to be nonlinear. Indeed, DP solutions only exist for the cases where the is an affine function of the MWPS. For example, if the exact penalty function \eqref{eq:zeta} is required, the DP scheme presented in Proposition \ref{prop:ClassicalDP} fails to apply.

These observations, however, do not exclude the existence of a solution to problem \eqref{P:original} via DP, but they exclude the existence of classical cost-to-go functions expressed in terms of the state $\vect s_k$ alone.
We now discuss how one can in principle deal with situations where Assumption \ref{ass:commutes} is violated and strict MWCC requirement should be satisfied. 
A solution via DP arguably requires one to augment the state space to enlarge the information at time $k$ involved in the decision-making.

To that end, let us define the sequence of functions $F_k:\mathcal{S}\to[0,\infty),$ for $k=0,1\ldots,N$ as follows:
\begin{align*}
	&F_0 (\vect s) = \mathbf{1}_{\mathbb S}(\vect s)\,,\\
	&F_k (\vect s)=\mathbb{P}\left[\,\vect s_{1,\ldots,k-1}\in\mathbb S \wedge \vect s_k=\vect{s}\,\middle|\,\vect s_0, \vect\pi\,\right],
\end{align*}
This sequence has the forward linear dynamics:
\begin{align}\label{eq:sysF}
	F_{k+1}(\vect s) = \int_{\mathbb S}F_{k}(\vect s')\rho\left[\vect s_{k+1}=\vect s\,|\,\vect s_k=\vect s',\vect a_k\right]\mathrm d \vect s'.
\end{align} 
Given a state $\vect s_0$, the associated MWPS then reads as $$\int_{\mathbb S} F_N(\vect s)\mathrm d\vect s=\mathbb{P}\left[\, \vect s_{1,\ldots,N}\in\mathbb S \,\middle|\,\vect s_0,\vect\pi \right].$$
At this point, by integrating the MWCC into the cost function and using Proposition \ref{pro:equiv}, we can then reformulate the MWCC-OCP problem \eqref{P:original} into its equivalent form
	\begin{align} \label{FinalProblem}
		J^{*}(\vect s_0)=\min_{\vect\pi}&\,\,\, \mathbb E \left[ \sum_{k=0}^{N-1}\ell(\vect s_k,\vect a_k) + \ell_N(\vect s_N)\,\middle|\,\vect s_0,\vect\pi  \right]\nonumber\\ 
		&\hspace*{0.4cm}+\zeta\left( \int_{\mathbb S} F_N(\vect s)\mathrm d\vect s\right)\,,
	\end{align}
where the function $\zeta$ is given by e.g. \eqref{eq:zeta}.

It is useful to observe that $\vect s_k$ is a vector in $\mathbb R^{n_s}$ stochastically decided by $\vect s_k$ and $\vect a_k$ through \eqref{eq:PDF}, while $F_k$ is a functional in some functional space deterministically decided by given $F_0$ and policy $\vect{\pi}_0,\dots,\vect\pi_{k-1}$ through \eqref{eq:sysF}. Let us denote  the augmented state as $\delta_k=\left(\vect s_k, F_k\right),$ consisting of a vector state and a functional state. The dynamics of state $\vect s_k$ and $F_k$ are given by \eqref{eq:PDF} and \eqref{eq:sysF}, respectively. 
Naturally, the control $\vect a_k$ should now depend on the new state $\delta_k$, or equivalently, a policy sequence should consist of policies $\vect\pi_k$ based on the regular state $\vect s_k$, as well as the functional state $F_k$.

By using the new augmented state, \eqref{FinalProblem} is identical to the basic problem format of finite-horizon optimal control problem, but includes a terminal constraint on the additional state $F_N$. Consequently, we can retain a DP solution to \eqref{FinalProblem} without the need of Assumption \ref{ass:commutes}, at the expense of a significantly more complex state-space to manipulate. The following Proposition details that DP solution.
\begin{proposition}\label{prop:OurDP}
	For every feasible initial state $\vect s_0$, the optimal cost $J^*(\vect s_0)$ of problem \eqref{FinalProblem} is equal to $J_0(\delta_0)$, given by the last step of the following backward-iteration algorithm:
	\begin{align}
		&J_N(\delta_N) = \ell_N(\vect s_N)+\zeta\left(\int_{\mathbb S} F_N(\vect s)\mathrm d\vect s\right)\,,\nonumber\\
		&J_k(\delta_k) = \min_{\vect a_k}\,\ell(\vect s_k,\vect a_k)\nonumber\\
		&\hspace{2.1cm}+\int_{\mathbb S} J_{k+1}(\delta_{k+1})\rho[\vect s_{k+1}|\vect s_k,\vect a_k]\mathrm{d} \vect s_{k+1}\,,\nonumber\\
		&\hspace{2.3cm}k=N-1,N-2,\ldots,0.\label{eq:OurDP_2}
	\end{align}
	Furthermore, if $\vect\pi_k^*(\delta_k)=\vect a_k^*$ minimizes the right hand side of \eqref{eq:OurDP_2} for each $\delta_k=(\vect s_k, F_k)$ and $k$, then the policy $\vect\pi^*=\{\vect\pi^*_0,\ldots,\vect\pi^*_{N-1}\}$ is optimal for problem \eqref{FinalProblem}. 
	\end{proposition}
	
\BeginProof 
For $k=N-1,N-2,\ldots,0$, let $J^*_k(\delta_k)$ be the optimal cost for the $(N-k)$-stage problem that starts at state $\delta_k=(\vect s_k, F_k)$ and time $k$, and ends at time $N$, that is, 
\begin{align*}
	J^*_k(\delta_k) = \min_{\vect \pi^k}\,\mathbb E \left[ \sum_{i=k}^{N-1}\ell(\vect s_i,\vect a_i) + \ell_N(\vect s_N)\,\middle|\,\vect s_k,\vect\pi^k  \right]
\end{align*}  
For $k = N$, we define 
$J^*_N(\delta_N) = \ell_N(\vect s_N)+\zeta\left( \int_{\mathbb S} F_N(\vect s)\mathrm d\vect s\right).$
We will show by induction that the functions $J^*_k$ are equal to the functions $J_k$ generated by the DP algorithm \eqref{eq:OurDP_2}, such that at $k=0$ the desired result will be obtained. 

By definition, we have $J^*_N = J_N = \ell_N+\zeta.$ Assume that for some $k+1$ and all $\delta_{k+1}$, we have $J^*_{k+1}(\delta_{k+1})=J_{k+1}(\delta_{k+1})$. Then, since $\vect\pi^k=(\vect\pi_k,\vect\pi^{k+1})$, we have for all $\delta_k$
\begin{align*}
	&J^*_k(\delta_k) =\min_{\vect \pi^k}\,\mathbb E \left[ \sum_{i=k}^{N-1}\ell(\vect s_i,\vect a_i) + \ell_N(\vect s_N)\,\middle|\,\vect s_k,\vect\pi^k  \right]\\
	&=\min_{\vect\pi_k}\,\ell(\vect s_k,\vect \pi_k(\delta_k))\\
	&+\mathbb{E}_{\vect s_{k+1}}\hspace{-0.1cm}\left[\min_{\vect\pi^{k+1}}\mathbb E \hspace{-0.05cm}\left[ \sum_{i=k+1}^{N-1}\ell(\vect s_i,\vect a_i) + \ell_N(\vect s_N)\,\middle|\,\vect s_{k+1},\vect\pi^{k+1}  \right]\right]
\end{align*}

Finally, by using the definition of $J^*_{k+1}(\delta_{k+1})$, we can further have
\begin{align*}
	&J^*_k(\delta_k)= \min_{\vect\pi_k}\,\ell(\vect s_k,\vect \pi_k(\delta_k)) + \mathbb E_{\vect s_{k+1}}\left[J^*_{k+1}(\delta_{k+1})  \right] \\	
	&= \min_{\vect\pi_k}\,\ell(\vect s_k,\vect \pi_k(\delta_k)) + \mathbb E_{\vect s_{k+1}}\left[J_{k+1}(\delta_{k+1})  \right] \\
	&= \min_{\vect a_k\in\mathcal A}\,\ell(\vect s_k,\vect a_k) + \int_{\mathbb S}J_{k+1}(\delta_{k+1}) \rho\left[\vect s_{k+1}\,\middle|\,\vect s_k, \vect a_k\right]\mathrm{d}\vect s_{k+1}\\
	&= J_k(\delta_k)\,,
\end{align*}
The explanations to the equations above are analogous to the corresponding ones given in the proof of Proposition \ref{prop:ClassicalDP}. Moreover, since $F_0(\vect s_0) \equiv 1$ for every feasible initial state $\vect s_0\in\mathbb S$, we have $J^*(\vect s_0) = J_0(\delta)$. 
\EndProof

\begin{remark}
	The DP principle given by Proposition \ref{prop:OurDP} actually holds for any penalty function $\zeta$, i.e., not just for the case of \eqref{eq:zeta}, which makes Proposition \ref{prop:OurDP} more universal. However, in order to enforce the exact MWCC, then $\zeta$ defined in \eqref{eq:zeta} can be used.
\end{remark}
 
As is typically the case, state augmentation often comes at a price of making very complex state and/or control spaces. This state augmentation has a complex state space that comprise a regular Euclidean space and a functional space. Nonetheless, the exact dynamic programming scheme proposed in Proposition \ref{prop:OurDP}, for the first time, gives the exact global optimal solution for the MWCC-OCP. Therefore, Proposition \ref{prop:OurDP} can serve as a stepping-stone for possibly many approximate dynamic programming methods to be developed in the future.

\section{Case Study}\label{sec:casestudy}
In this section, we analytically demonstrate how to use the DP algorithm proposed in the preceding section on a simple example. 
Consider the one-dimensional system of the form
\begin{align*}
	\vect s_{k+1} = \vect s_k + \vect a_k + \vect w_k,\qquad k=0,1,
\end{align*}
with initial state $\vect s_0 =\bar{\vect s}$ which is supposed to be feasible regarding the following problem setup.
Here, $\vect a_k\in[-0.1,0.1]$ and disturbance $\vect w_k\sim\mathcal N(0,0.0001)$. The safe set $\mathbb S=[-1,1]$. The mission duration $N=2$. The cost functions $\ell(\vect{s},\vect{a})= \vect{s}^2+\vect{a}^2$ and $\ell_2(\vect{s})= \vect{s}^2$. The exact penalty function $\zeta$ follows from \eqref{eq:zeta} with risk bound $\varepsilon=0.1$. We observe that the augmented functional states $F_1,\, F_2,\, F_3$ are given by: 
\begin{align*}
	&F_0 (\vect s_0) = \mathbf{1}_{[-1,1]}(\bar{\vect s}) = 1\,,\\
	&F_1 (\vect s)=F_0 (\bar{\vect s})\rho\left[ \vect s_1=\vect{s}\,\middle|\,\vect s_0=\bar{\vect s}, \vect a_0=\vect{ x_0}\,\right] \\
	&\hspace*{0.8cm}= \frac{1}{0.01\sqrt{2\pi}}e^{-\frac{(\vect s-\bar{\vect s}-\vect x_0)^2}{0.002}}\\
	&F_2 (\vect s)= \int_{-1}^{1}F_1(\vect s')\rho\left[ \vect s_2=\vect{s}\,\middle|\,\vect s_1=\vect s', \vect a_1=\vect{ x}_1\,\right] \mathrm{d}\vect s'\\
	&\hspace*{0.82cm}=\int_{-1}^{1} \frac{1}{(0.01\sqrt{2\pi})^2}e^{-\frac{(\vect s'-\bar{\vect s}-\vect x_0)^2+(\vect s-\vect s'-\vect x_1)^2}{0.002}}\mathrm{d}\vect s'
\end{align*}
Note that the functions $F_1,\, F_2$ are functions involving parameters $\vect x_0,\,\vect x_1 \in[-0.1,0.1]$.  


By using the augmented state $\delta_{k} = (\vect s_k, F_k)$, the DP algorithm takes the following form. 
At stage $2$, we initialize the value function
\begin{align*}
	J_2(\delta_2) &= J_2(\vect s_2, F_2) = \ell_N(\vect s_2)+\zeta\left(\int_{-1}^1 F_2(\vect s)\mathrm d\vect s\right)\\
	&= 
	\begin{cases}
		\vect s^2_2,        &\text{if}\,\, \zeta\left(\int_{-1}^1 F_2(\vect s)\mathrm d\vect s\right)\geq 0.9 \\
		\infty,   &\text{otherwise}\,.
	\end{cases}
\end{align*}
for all possible state $\delta_2$, At stage $1$, we solve the following optimization problem for all possible $\delta_{1}$ 
\begin{align*}
	&J_1(\delta_1) = J_1(\vect s_1, F_1)\\
	&=\min_{\vect a_1\in[-0.1,0.1]}\,\ell(\vect s_1,\vect a_1)+ \mathbb E_{\vect s_2} \left[ J_2(\delta_2)\mid \vect s_1, \vect a_1 \right] \\
	&=\min_{\vect a_1\in[-0.1,0.1]}\,\vect s_1^2+\vect a_1^2\\
	& \hspace*{1.3cm}+\int_{-1}^{1}\frac{1}{0.01\sqrt{2\pi}}e^{-\frac{(\vect s_2-\vect s_1-\vect a_1)^2}{0.002}}\cdot J_2(\delta_2) \mathrm{d}\vect s_2\,, 
\end{align*}
for all possible state $\delta_1$.
The optimal control policy $\vect \pi_1$ is given by solving the above optimization problem.

At stage $0$, for the fixed initial state $\vect s_0=\bar{\vect s}$, the optimal control input $\vect \pi_0(\vect s_0)$ is given by solving the following optimization problem  
\begin{align*}
	&J^*(\vect s_0) = J_0(\delta_0) = J_0(\bar{\vect s}, F_0)\\
	&=\min_{\vect a_0\in[-0.1,0.1]}\,\ell(\vect s_0,\vect a_0)+ \mathbb E_{\vect s_1} \left[ J_1(\delta_1)\mid \vect s_0=\bar{\vect s}, \vect a_0 \right]\\
	&=\min_{\vect a_0\in[-0.1,0.1]}\,\vect a_0^2+\int_{-1}^1 J_{1}(\delta_{1}) \frac{1}{0.01\sqrt{2\pi}}e^{-\frac{(\vect s_1-\bar{\vect s}-\vect a_0)^2}{0.002}} \mathrm{d} \vect s_{1}
\end{align*}
 
This ``toy" example illustrates the essential procedures of using Proposition \ref{prop:OurDP}. A practical implementation on a more complete example is beyond the scope of this paper and is being investigated in our current work.

\section{CONCLUSIONS}\label{sec:conclusion}
In this paper, we investigate solutions for mission-wide chance-constrained optimal control problems via Dynamic Programming. We show that classic Dynamic Programming recursions on the state-space of the problem are possible if the penalty imposed on mission-wide chance constraints violations commutes with the expected value operator underlying the stochastic dynamics. We show that this requirement is not fulfilled when imposing hard mission-wide chance constraints. We then present a state augmentation that tackles the problem. The resulting augmented space consists of a regular Euclidean space and a functional space. The proposed dynamic programming scheme for the mission-wide chance-constrained optimal control problems can hopefully play a fundamental role for developing approximation methods, because it characterizes the optimal solutions.


%




%
%
%

\bibliographystyle{ieeetr}

\bibliography{reference_dp.bib}

\end{document}